# Can Turing machines capture everything we can compute?

Bhupinder Singh Anand

If we define classical foundational concepts constructively, and introduce non-algorithmic effective methods into classical mathematics, then we can bridge the chasm between truth and provability, and define computational methods that are not Turing computable.

## Contents

## 1. Introduction

In a short opinion paper, "Computation Beyond Turing Machines" [WG03], Peter Wegner and Dina Goldin advance the thesis that:

> "A paradigm shift is necessary in our notion of computational problem solving, so it can provide a complete model for the services of today's computing systems and software agents."

We note that Wegner and Goldin's arguments, in support of their thesis, seem to reflect an extraordinarily eclectic view of mathematics, combining both acceptance of, and frustration at, the standard interpretations and dogmas of classical mathematical theory:

(*i*) ...Turing machines are inappropriate as a universal foundation for computational problem solving, and ... computer science is a fundamentally non-mathematical discipline.

(*ii*) (*Turing's*) 1936 paper ... proved that mathematics could not be completely modeled by computers.

(*iii*) ...the Church-Turing Thesis ... equated logic, lambda calculus, Turing machines, and algorithmic computing as equivalent mechanisms of problem solving.

(*iv*) Turing implied in his 1936 paper that Turing machines ... could not provide a model for all forms of mathematics.

(*v*) ...Gödel had shown in 1931 that logic cannot model mathematics ... and Turing showed that neither logic nor algorithms can completely model computing and human thought.

These remarks vividly illustrate the Catch-22 situation, noted in Anand [An03c], with which not only Theoretical Computer Sciences, but all applied sciences that depend on mathematics - for providing a verifiable language to express their observations precisely - are faced:

Are formal classical theories essentially unable to adequately express the extent and range of human cognition, or does the problem lie in the way formal theories are classically interpreted at the moment? The former addresses the question of whether there are absolute limits on our capacity to express human cognition unambiguously; the latter, whether there are only temporal limits - not necessarily absolute - to the capacity of classical interpretations to communicate unambiguously that which we intended to capture within our formal expression.

Prima facie, applied science continues, perforce, to interpret mathematical concepts Platonically, whilst waiting for mathematics to provide suitable, and hopefully reliable, answers as to how best it may faithfully express its observations verifiably.

This dilemma is also reflected in Fortnow's on-line rebuttal[1] of Wegner and Goldin's thesis, and of their reasoning. Thus Fortnow divides his faith between the standard interpretations of classical mathematics (and, possibly, the standard set-theoretical models of formal systems such as standard Peano Arithmetic), and the classical computational theory of Turing machines. He relies on the former to provide all the proofs that matter:

> "Not every mathematical statement has a logical proof, but logic does capture everything we can prove in mathematics, which is really what matters";

and, on the latter to take care of all essential, non-provable, truth:

> "...what we can compute is what computer science is all about".

However, as we argue later, Fortnow's faith in a classical Church-Turing Thesis that ensures:

> "...Turing machines capture everything we can compute",

may be as misplaced as his faith in the infallibility of standard interpretations of classical mathematics.

The reason: There are, prima facie, reasonably strong arguments for a Kuhnian paradigm shift; not, as Wegner and Goldin believe, in the notion of computational problem solving, but in the standard interpretations of classical mathematical concepts.

---

[1] Cf. Lance Fortnow's Web Log <http://fortnow.com/lance/complog/archive/2003_04_06_archive.html>.

However, Wegner and Goldin could be right in arguing that the direction of such a shift must be towards the incorporation of non-algorithmic effective methods into classical mathematical theory; presuming, from the following remarks, that this is, indeed, what "external interactions" are assumed to provide beyond classical Turing-computability:

> (*vi*) ...that Turing machine models could completely describe all forms of computation ... contradicted Turing's assertion that Turing machines could only formalize algorithmic problem solving ... and became a dogmatic principle of the theory of computation.

> (*vii*) ...interaction between the program and the world (environment) that takes place during the computation plays a key role that cannot be replaced by any set of inputs determined prior to the computation.

> (*viii*) ...a theory of concurrency and interaction requires a new conceptual framework, not just a refinement of what we find natural for sequential [algorithmic] computing.

> (*ix*) ...the assumption that all of computation can be algorithmically specified is still widely accepted.

A widespread notion of particular interest, which seems to be recurrently implicit in Wegner and Goldin's assertions too, is that mathematics is a dispensable tool of science, rather than its indispensable mother tongue. In Anand [An02] we argue that the roots of such beliefs may lie in ambiguities, in the classical definitions of foundational elements, that allow the introduction of non-constructive - hence non-verifiable, non-computational, ambiguous, and essentially Platonic - elements into the standard interpretations of classical mathematics[2].

---

[2] See, for instance, the arguments in Anand [An03b].

Consequently, standard interpretations of classical theory may, inadvertently, be weakening a desirable perception - of mathematics as the lingua franca of scientific expression - by ignoring the possibility that, since mathematics is, indeed, indisputably accepted as the language that most effectively expresses and communicates intuitive truth, the chasm between formal truth and provability must, of necessity, be bridgeable.[3]

We now consider one such bridge.

## 2. Non-constructivity in classical theory

We note that, in [Me90], Mendelson's following remarks (*italicised parenthetical qualifications added*), implicitly imply that classical definitions of various foundational elements can be argued as being either ambiguous, or non-constructive, or both:

> "Here is the main conclusion I wish to draw: it is completely unwarranted to say that CT (*Church's Thesis*) is unprovable just because it states an equivalence between a vague, imprecise notion (effectively computable function) and a precise mathematical notion (partial-recursive function). ... The concepts and assumptions that support the notion of partial-recursive function are, in an essential way, no less vague and imprecise (*non-constructive, and intuitionistically objectionable*) than the notion of effectively computable function; the former are just more familiar and are part of a respectable theory with connections to other parts of logic and mathematics. (The notion of effectively computable function could have been incorporated into an axiomatic presentation of classical mathematics, but the acceptance of CT made this unnecessary.) ... Functions are defined in terms of sets, but the concept of set is no clearer (*not more non-constructive, and intuitionistically objectionable*) than that of function and a foundation of mathematics can be based on a theory using function as primitive notion instead of set. Tarski's definition of truth is formulated in set-theoretic terms, but the notion of set is no clearer (*not more non-constructive, and intuitionistically objectionable*), than that of truth. The model-theoretic definition of logical validity is based ultimately on set theory, the foundations of which are no clearer (*not more non-constructive, and intuitionistically objectionable*) than our intuitive (*non-constructive, and intuitionistically objectionable*) understanding of logical validity. ... The notion of Turing-computable function is no clearer (*not more non-constructive, and intuitionistically objectionable*) than, nor more mathematically useful (foundationally speaking) than, the notion of an effectively computable function."[4]

---

[3] Of interest is Davis' argument [Da95] that an unprovable truth may, indeed, be arrived at algorithmically.

## 3. Constructive definitions

We further note that, if we accept that there can be non-algorithmic effective methods that are applicable individually, but not uniformly (algorithmically), then we can constructively express some of these classical notions[5] of various foundational elements such as "mathematical object", "effective computability", "truth of a formula under an interpretation", "set", "Church's Thesis" etc., in terms of a smaller number of primitive - formally undefined but intuitively well-accepted as unambiguous - mathematical terms as below[6]:

(*i*) **Primitive mathematical object**: A primitive mathematical object is any symbol for an individual constant, predicate letter, or a function letter, which is defined as a primitive symbol of a formal mathematical language.

---

[4] Reproduced from Bringsjord [Br93].

[5] We take Mendelson [Me64] as a standard exposition of classical theory and its standard interpretations.

[6] All the technical terms are defined explicitly in Anand [An02].

(*ii*) **Formal mathematical object**: A formal mathematical object is any symbol for an individual constant, predicate letter, or a function letter that is either a primitive mathematical object, or that can be introduced through definition into a formal mathematical language without inviting inconsistency.

(*iii*) **Mathematical object**: A mathematical object is any symbol that is either a primitive mathematical object, or a formal mathematical object.

(*iv*) **Set**: A set is the range of any function whose function letter is a mathematical object.

(*v*) **Individual computability**: A number-theoretic function $F(x)$ is individually computable if, and only if, given any natural number $k$, there is an individually effective method (which may depend on the value $k$) to compute $F(k)$.

(*vi*) **Uniform computability**: A number-theoretic function $F(x)$ is uniformly computable if, and only if, there is a uniformly effective method (necessarily independent of $x$) such that, given any natural number $k$, it can compute $F(k)$.

(*vii*) **Effective computability**: A number-theoretic function is effectively computable if, and only if, it is either individually computable, or it is uniformly computable.

(*viii*) **Individual truth**: A string $[F(x)]$[7] of a formal system P is individually true under an interpretation M of P if, and only if, given any value $k$ in M, there is an individually effective method (which may depend on the value $k$) to determine that the interpreted proposition $F(k)$ is satisfied in M ([Me64], p49-52).

---

[7] We use square brackets to distinguish between the uninterpreted string [F] of a formal system, and the symbolic expression "F" that corresponds to it under a given interpretation that unambiguously assigns formal, or intuitive, meanings to each individual symbol of the expression "F".

(*ix*) **Uniform truth**: A string $[F(x)]$ of a formal system P is uniformly true under an interpretation M of P if, and only if, there is a uniformly effective method (necessarily independent of $x$) such that, given any value $k$ in M, it can determine that the interpreted proposition $F(k)$ is satisfied in M.

(*x*) **Effective truth**: A string $[F(x)]$ of a formal system P is effectively true under an interpretation M of P if, and only if, it is either individually true in M, or it is uniformly true in M.

(*xi*) **Individual Church Thesis**: If, for a given relation $R(x)$, and any element $k$ in some interpretation M of a formal system P, there is an individually effective method such that it will determine whether $R(k)$ holds in M or not, then every element of the domain D of M is the interpretation of some term of P, and there is some P-formula $[R'(x)]$ such that:

$R(k)$ holds in M if, and only if, $[R'(k)]$ is P-provable.

(In other words, the Individual Church Thesis postulates that, if a relation $R$ is effectively decidable individually (possibly non-algorithmically) in an interpretation M of some formal system P, then $R$ is expressible in P, and its domain necessarily consists of only mathematical objects, even if the predicate letter $R$ is not, itself, a mathematical object.)

(*xii*) **Uniform Church Thesis**: If, in some interpretation M of a formal system P, there is a uniformly effective method such that, for a given relation $R(x)$, and any element $k$ in M, it will determine whether $R(k)$ holds in M or not, then $R(x)$ is the interpretation in M of a P-formula $[R(x)]$, and:

$R(k)$ holds in M if, and only if, $[R(k)]$ is P-provable.

(Thus, the Uniform Church Thesis postulates that, if a relation *R* is effectively decidable uniformly (necessarily algorithmically) in an interpretation M of a formal system P, then, firstly, *R* is expressible in P, and, secondly, the predicate letter *R*, and all the elements in the domain of the relation *R*, are necessarily mathematical objects.)

## 4. Effective solvability of the Halting problem

The significance of defining classical foundational concepts in a constructive, and, prima facie, in an intuitionistically unobjectionable way, is seen in the following argument.

**Theorem 1**: The Uniform Church Thesis implies that the Halting problem is effectively solvable.

**Proof**: If a number-theoretic function *F*(*x*) is Turing-computable, then it is partial recursive ([Me64], p233, Corollary 5.13). We may thus assume that such an *F* is obtained from a recursive function *G* by means of the unrestricted $\mu$-operator; in other words, that (cf. [Me64], p214):

$F(x) = \mu y(G(x, y) = 0)$.

If [*H*(*x*, *y*)] expresses ~(*G*(*x*, *y*) = 0) in a formal system of Arithmetic such as standard PA, we then consider the PA-provability, and truth in the standard interpretation M of PA, of the formula [*H*(*a*, *y*)] for a given numeral [*a*] of PA, as below:

(*a*) Let Q1 be the meta-assertion that [*H*(*a*, *y*)] is not effectively true in M. Hence there is no effective method in M to determine that, for any given *y* in M, *y* satisfies [*H*(*a*, *y*)] classically. It follows that there is no uniformly effective method (algorithm/Turing machine) in M to determine that, for any given *y* in M, *y* satisfies [*H*(*a*, *y*)] classically.

Since *G*(*a*, *y*) is recursive, it follows that there is some finite *k* such that any Turing machine T1(*y*) that computes *G*(*a*, *y*) will halt and return the value 0 for *y* = *k*.

(*b*) Next, let Q2 be the meta-assertion that $[H(a, y)]$ is effectively true in the standard interpretation M of PA, but that there is no uniformly effective method (algorithm/Turing machine) in M to determine that, for any given $y$ in M, $y$ satisfies $[H(a, y)]$ classically.

Since $G(a, y)$ is recursive, it follows that there is some finite $k$ such that the T1($y$) will halt, and return a symbol for self-termination (looping[8]) for $y = k$.

(*c*) Finally, let Q3 be the meta-assertion that $[H(a, y)]$ is effectively true in the standard interpretation M of PA, and that there is a uniformly effective method in M to determine that, for any given $y$ in M, $y$ satisfies $[H(a, y)]$ classically. We then have that that $[H(a, y)]$ is uniformly true in the standard interpretation M of PA.

If we assume a Uniform Church Thesis, it then follows ([An02], Meta-lemma 8) that $[H(a, y)]$ is PA-provable. Let $h$ be the Gödel-number of $[H(a, y)]$. We consider, now, Gödel's primitive recursive number-theoretic relation $x$B$y$ ([Go31a], p22, definition 45), which holds in M if, and only if, $x$ is the Gödel-number of a proof sequence in PA for the PA-formula whose Gödel-number is $y$. It follows that there is some finite $k$ such that any Turing machine T2($y$), which computes the characteristic function of $x$B$h$, will halt and return the value 0 for $x = k$.

Since Q1, Q2, and Q3 are mutually exclusive and exhaustive, it follows that, when run simultaneously over the sequence 1, 2, 3, ... of values for $y$, one of {T1(y) // T2(y)} will always halt for some finite value of $y$ for any given $a$.

Thus, the Halting problem is effectively solvable if we assume a Uniform Church Thesis.¶

---

[8] We note that any Turing machine can be designed to recognise a "looping" situation; it simply records every instantaneous tape description at the execution of each machine instruction, and compares the current instantaneous tape description with the record. It can thus be meta-programmed to abort a loop, and return a meta-symbol indicating self-termination.

**Corollary 1.1**: The Uniform Church Thesis implies that the parallel duo of Turing machines {T1(y) // T2(y)} is not a Turing machine.

**Corollary 1.2**: The Uniform Church Thesis implies that the classical Church-Turing thesis is false.

## 5. Is there a case for a Uniform Church Thesis?

The question of whether, as Fortnow believes, Turing machines can "capture everything we can compute", or whether, as Wegner and Goldin suggest, they are "inappropriate as a universal foundation for computational problem solving", could, thus, be answered decisively if there were, indeed, a case for introducing non-algorithmic effective methods into classical mathematics, and for introducing constructive foundational concepts such as that of a Uniform Church Thesis. We address this issue in Anand [An02], Anand [An03a] and Anand [An03c].

(*Updated: Saturday 10th May 2003 7:44:10 AM IST by re@alixcomsi.com*)